\documentclass{amsart}
\newtheorem{thm}{Theorem}[section]
\newtheorem{cor}[thm]{Corollary}
\newtheorem{lem}[thm]{Lemma}
\newtheorem{prop}[thm]{Proposition}
\theoremstyle{definition}
\newtheorem{defn}[thm]{Definition}
\theoremstyle{remark}

\newtheorem{exmp}{Example}
\numberwithin{equation}{section}
\newcommand{\scalar}[2]{\langle#1,#2\rangle}
\newcommand{\Real}{\mathbb{R}}

\newcommand{\To}{\longrightarrow}

\newcommand{\pars}[1]{\left(#1\right)}
\newcommand{\brck}[2]{\{#1,#2\}}
\newcommand{\set}[1]{\left\{#1\right\}}
\newcommand{\jami}[2]{\sum\limits_{#1}^{#2}}
\newcommand{\all}[1]{\forall\,#1}

\newcommand{\gamomd}{\,\Rightarrow\,}

\newcommand{\smooth}[1]{\ensuremath{C^\infty(#1)}}
\newcommand{\fields}[2]{\mathcal{V}^{#1}(#2)}
\newcommand{\vecf}[1]{\fields{1}{#1}}
\newcommand{\forms}[2]{\ensuremath{\Omega^{#1}(#2)}}
\newcommand{\allforms}[1]{\Omega(#1)}

\newcommand{\tolfas}{\,\Leftrightarrow\,}
\newcommand{\lder}[2]{L_{#1}(#2)}
\newcommand{\talg}[1]{\widetilde{#1}}
\newcommand{\distr}[1]{\ensuremath{\operatorname{D}'(#1)}}
\newcommand{\distro}[1]{\ensuremath{\operatorname{D}_0'(#1)}}
\newcommand{\distrf}[1]{\ensuremath{\Omega'(#1)}}
\newcommand{\distrfo}[1]{\ensuremath{\Omega_0'(#1)}}
\newcommand{\genc}[1]{\ensuremath{\operatorname{D}_0'(#1)}}
\newcommand{\emb}[2]{\ensuremath{\mathcal{E}_{_{#1}}(#2)}}
\newcommand{\embm}[1]{\ensuremath{\mathcal{E}_{_{#1}}}}
\begin{document}

\title{Bott Connection and Generalized Functions on Poisson Manifold}%
\author{Zakaria Giunashvili}%
\address{Department of Theoretical Physics,
         Institute of Mathematics,
         Georgian Academy of Sciences,
         Tbilisi, Georgia}%
\email{zaqro@gtu.edu.ge}%

\thanks{This work was supported by Nishina Memorial Foundation and is based on the talks given by the author
        at Waseda University and Keio University (Tokyo, Japan). The author gratefully acknowledges the
        hospitality of Yokohama City University where the work was completed.}%

\date{January 31, 2003}%
% ----------------------------------------------------------------
\begin{abstract}
We extend the problem of finding Hamiltonian-invariant volume forms on a Poisson manifold to the problem of
construction of Hamiltonian-invariant generalized functions. For this we introduce the notion of generalized
center of a Poisson algebra, which is the space of generalized Casimir functions. We study as the case
when the set of test-objects for generalized functions is the space of compactly supported smooth functions, so the case
when the test-objects are $n$-forms, where $n$ is the dimension of the Poisson manifold. We describe the
relations of this problem with the homological properties of the Poisson structure, with Bott connection for the
corresponding symplectic foliation and the modular class.
\end{abstract}
\maketitle
% ----------------------------------------------------------------
%%%%%%%%%%%%%%%%%%%%%%%%%%%%%%%%%%%%%%%%%%%%%%%%%%%
\section{Introduction}
If the rank of a Poisson structure on a smooth manifold $M$ is constant, there are various tools in the theory of
foliations that can be applied to study the corresponding symplectic foliation. But there are many even
very simple degenerate Poisson structures that are singular, i.e, the rank of the corresponding contravariant
tensor field varies. Therefore, some new methods are needed to handle this situation. We define an action
of the Poisson algebra $\smooth{M}$ on the space of generalized functions (distributions). The
main idea of our construction is the following: if the center of a Poisson algebra is more than
the constants, then the Poisson structure is certainly degenerate, but the converse is not true,
in general. We introduce the notion of generalized center, which is the set of such distributions
on a Poisson manifold that commute with every ordinary function. We show that the support of the
generalized central elements (generalized Casimir functions) consists of the symplectic leaves of
the Poisson manifold. We show that every compact symplectic leaf gives a family of generalized
Casimirs: a basic one and its derivatives by some special vector fields which are transversal to
the symplectic leaf and flat with respect to the Bott connection.
%%%%%%%%%%%%%%%%%%%%%%%%%%%%%%%%%%%%%%%%%%%%%%%%%%%
\section{Algebraic Definition of Bott Connection}
\newcommand{\der}{Der(A)}
\newcommand{\deri}{Der_I(A)}
\newcommand{\derio}{Der_I(A)_0}
\newcommand{\gami}{\Gamma(I)}
\newcommand{\gamio}{\Gamma(I)_0}
\newcommand{\inter}{\cap}
%%%%%%%%%%%%%%%%%%%%%%%%%%%%%%%%
Let $A$ be a commutative algebra and $Der(A)$ be the Lie algebra of derivations of $A$. As the algebra $A$ is
commutative the space $Der(A)$ is a module over $A$. For a given ideal $I$ in $A$ let us define the following
Lie subalgebras and  $A$-submodules in $Der(A)$
$$
\deri=\set{X\in\der\ |\ X(I)\subset I}
$$
and its subalgebra
$$
\derio=\set{Y\in\der\ |\ Y(A)\subset I}
$$
There is a canonical homomorphism
$$
\rho:\deri\To Der(A/I)
$$
and as it follows from the definition, we have that $\derio=\ker(\rho)$. When $\rho$ is an epimorphism the ideal $I$ is
said to be a \emph{submanifold constraint ideal} (see \cite{Masson}). For brevity we shall use the term
\emph{submanifol ideal}. In this case we have the isomorphism
$$
Der(A/I)\cong\frac{\deri}{\derio}
$$
We use the term \emph{distribution} for any Lie subalgebra which is also an $A$-submodule of $\der$.
\begin{defn}[Integral of Distribution]
An ideal $I$ in $A$ is said to be an \emph{integral} for a given distribution $D$ if $D\subset\deri$.
\end{defn}
Let us denote the intersection $D\inter\derio$ by $D_0$.
\begin{defn}[Complete Integral]
An ideal $I\subset A$ is said to be a \emph{complete integral} for a given distribution $D\subset\der$ if the
following three conditions are satisfied
\begin{itemize}
\item[(1)] $I$ is a submanifold ideal for $A$
\item[(2)] $I$ is an integral for $D$
\item[(3)] $D/D_0=\frac{\deri}{\derio}\quad$ (\ $\cong Der(A/I)$ because of (1))
\end{itemize}
\end{defn}
Let $I\cdot D=\set{\sum a_iX_i\ |\ a_i\in I,\ X_i\in D}$. It is clear that $I\cdot D$ is a subspace of $D_0$.
\begin{defn}[Regular Integral]
A complete integral $I$ for a given distribution $D\in\der$ is said to be \emph{regular} if $D_0=I\cdot D$.
\end{defn}
The above definition of a regular integral is in agreement with the following theorem from the classical theory
of distributions
\begin{thm}[see \cite{Dazord}]
Let $D$ be a completely integrable distribution on a smooth manifold $M$, and $x_0\in M$ be a point where
$\dim(D_{x_0})=r$. Then the point $x_0$ has an open neighborhood $U$ with coordinates $(u_1,\ldots,u_n)$ such
that
$$
\set{\frac{\partial}{\partial u_i}(x)\ |\ i=1,\ldots,r}\subset D_x,\ \all x\in U
$$
\end{thm}
It follows from this theorem that if $x_0$ is a regular point for $D$, then $\dim(D_{x_0})$ is locally maximal,
and therefore, the set $\set{\frac{\partial}{\partial u_i}\ |\ i=1,\ldots,r}$ is a local basis of the
$\smooth{M}$-module $D$. If $N\subset M$ is the integral submanifold for $D$, passing through the point $x_0$
and $X\in D$ is such that $X|_N=0$, then
$X=\jami{i=1}{r}\varphi_i\frac{\partial}{\partial u_i},\ \varphi_i\in I_N$, where $I_N$ denotes the ideal of such
smooth functions on $M$ that vanish on $N$.

For a given ideal $I\in A$ let $\pi:A\To A/I$ be the natural quotient map. Introduce the following $A$-modules
$$
\gami=\set{X:A\To A/I\ |\ X(ab)=X(a)\pi(b)+\pi(a)X(b)}
$$
and its submodule
$$
\gamio=\set{Y\in\gami\ |\ Y(I)=\set{0}}
$$
We call the quotient module
$$
V(I)=\frac{\gami}{\gamio}
$$
the space of \emph{transversal derivations}. It is a module over the quotient algebra $A/I$. Consider the
following operation $[\ ,\ ]:\deri\times\gami\To\gami$, defined as
$$
[U,s]=\rho(U)\circ s-s\circ U,\quad U\in\deri,\ s\in\gami
$$
The submodule $\gamio$ is invariant for this operation and therefore we can reduce this operation to
$$
[\ ,\ ]:\deri\times V(I)\To V(I)
$$
The following properties of this operation follow directly from the definition
\begin{equation}\label{bracket_props}
\begin{array}{ll}
(1) & [aU,v]=a[U,v]\\
& \\
(2) & [U,kv]=U(k)v+k[U,v]
\end{array}
\end{equation}
where $a\in A,\ k\in A/I,\ U\in\deri$ and $v\in V(I)$.
\begin{defn}[Bott Connection]
\emph{Bott Connection} for $I\subset A$ which is a \emph{regular integral} for a distribution $D\in\der$, is an
action (covariant derivation) of the elements of Lie algebra $Der(A/I)\cong\frac{D}{I\cdot D}$ on the
$A/I$-module $V(I)$, defined as
$$
\nabla_X(v)=[\talg{X},v]
$$
where $X\in\frac{D}{I\cdot D},\ v\in V(I)$ and $\talg{X}\in D$ is an extension of $X$: $\rho(\talg{X})=X$.
\end{defn}
\textbf{Independence from the extension $\talg{X}$:} if $\talg{X}_1$ and $\talg{X}_2$ in $D$ are such that
$\rho(\talg{X}_1)=\rho(\talg{X}_2)$ then we have
$$
\begin{array}{c}
\talg{X}_1-\talg{X}_2\ \in I\cdot D\tolfas\talg{X}_1-\talg{X}_2=\sum k_iU_i,\ k_i\in I,\ U_i\in D\gamomd\\
\\
\gamomd[\talg{X}_1-\talg{X}_2,v]=[\sum k_iU_i,v]=\sum k_i[U_i,v]=0\ \in\ V(I)=\frac{\gami}{\gamio}
\end{array}
$$
\medskip\\
The properties of covariant derivation for $\nabla$ easily follow from the properies of the operation $[\ ,\ ]$.
In the case of a singular (i.e., not regular) integral, the definition of Bott connection encounters the following
problem: the subalgebra $D_0$, in general, does not coincide with $I\cdot D$, i.e., not every $\talg{X}\in D$ that
vanishes on $A/I$ is of the form $\talg{X}=\sum k_iU_i$ with $k_i\in I$ and $U_i\in D$. Hence the independence
of the expression $[\talg{X},v]$ from the extension $\talg{X}$ is problematic. A Poisson structure helps to
somehow resolve the difficulty. Some authors (see \cite{Ginzburg}) define the covariant derivation of a
transversal vector field not by tangent vectors but by 1-forms. In the next section we shall review this construction in
the case of a Poisson manifold.
%%%%%%%%%%%%%%%%%%%%%%%%%%%%%%%%%%%%%%%%%%%%%%%%%%%
\section{Bott Connection on a Symplectiv Leaf of Poisson Manifold}
\newcommand{\ortho}{T(N)^\perp}
\newcommand{\restr}{T^*(M)|_N}
%%%%%%%%%%%%%%%%%%%
Let $(M,\brck{\ }{\ })$ be a Poisson manifold. For any function $f\in\smooth{M}$, we denote by $X_f$ the
Hamiltonian vector field on $M$ corresponding to $f$. The antisymmetric contravariant tensor (bivector) field on
$M$ corresponding to the Poisson structure, defines a homomorphism of the vector bundles $T^*(M)\To T(M)$.
We denote the tangent vector corresponding to the cotangent vector $\alpha\in T^*(M)$ also by $X_\alpha$.
For a differential 1-form $\omega$ we have a vector field $X_\omega$ defined by the above homomorphism.

Let $N\subset M$ be a symplectic leaf. Denote by $\ortho$ the subbundle of the restricted bundle $\restr$:
$$
\ortho=\set{\omega\in T_x^*(M)\ |\ x\in N,\ \omega(T_x(N))=0}
$$
\emph{Bott connection} on the vector bundle $\ortho$ is a covariant derivation of the sections of this bundle by
the sections of the vector bundle $\restr$ defined as follows
$$
\nabla_\alpha(\omega)=\lder{X_{\talg{\alpha}}}{\talg{\omega}}|_N
$$
where:
\begin{itemize}
\item
	$\alpha$ is a section of $\restr$ and $\talg{\alpha}$ is its extension (at least local) on $M$:
	$\talg{\alpha}|_N=\alpha$;
\item
	$\omega$ is a section of $\ortho$ and $\talg{\omega}$ is its extension (at least local) on M:
	$\talg{\omega}|_N=\omega$;
\item
	$\lder{\cdot}{\cdot}$ is the operation of Lie derivation.
\end{itemize}
If $\talg{\alpha}=\varphi\cdot\talg{\alpha}_1$, where $\varphi\in\smooth{M}$ and $\talg{\alpha}_1$ is a 1-form,
then by definition of the Lie derivation we have
$$
\lder{X_{\talg{\alpha}}}{\talg{\omega}}|_N=(d\varphi\cdot\talg{\omega}(X_{\talg{\alpha}_1}))|_N+(\varphi\cdot\lder{X_{\talg{\alpha}_1}}{\talg{\omega}})|_N
$$
The first term of the RH side is 0 because the vector field $X_{\talg{\alpha}_1}|_N$ is always tangent to the
symplectic leaf $N$, and $\talg{\omega}|_N=\omega$ is a section of $\ortho$. Therefore we have
$$
\lder{X_{\talg{\alpha}_1}}{\talg{\omega}}|_N=\varphi|_N\cdot\lder{X_{\talg{\alpha}_1}}{\talg{\omega}}|_N
$$
from which easily follows as the independence of the expression $\lder{X_{\talg{\alpha}}}{\talg{\omega}}|_N$ from
the extension $\talg{\alpha}$ so the $\smooth{N}$-linearity of $\nabla_\alpha(\omega)$ by the argument $\alpha$.

\newcommand{\zxa}{X_{\talg{\alpha}}}
For $\talg{\omega}=\varphi\cdot\talg{\omega}_1$, we have
\begin{equation}\label{leibniz}
\lder{\zxa}{\varphi\cdot\talg{\omega}_1}|_N=(\zxa(\varphi)\cdot\talg{\omega}_1)|_N+(\varphi\cdot\lder{\zxa}{\talg{\omega}_1})|_N
\end{equation}
From this we obtain the independence from the extension $\talg{\omega}$ in the following way: if
$\talg{\omega}|_N=\talg{\omega}'|_N$, then $\talg{\omega}-\talg{\omega}'=\sum\varphi_i\talg{\omega}_i$, where
$\varphi_i|_N=0$; then recall that as $\zxa|_N$ is tangent to $N$, we have that $\zxa(\varphi_i)|_N=0$. From the
equality \ref{leibniz} also follows the Leibniz rule for the argument $\omega$ in $\nabla_\alpha(\omega)$.

Let us describe the dual definition of Bott connection, which is more compatible with the algebraic definition
introduced in the previous section. In this case, instead of the vector bundle $\ortho$, we take the quotient
bundle $\frac{T(M)|_N}{T(N)}\equiv V(N)$, which is the bundle of dual spaces of the fibers of $\ortho$. For a
section $\alpha$ of $\restr$ and a section $s$ of $V(N)$, we set
\begin{equation}\label{bott}
\nabla_\alpha(s)=\rho([\zxa,\talg{s}])
\end{equation}
where: $\talg{\alpha}$ is an extension of $\alpha$ as in the previous definition; $\talg{s}$ is a vector field on
$M$, such that $\rho(\talg{s}|_N)=s$ and $\rho:T(M)|_N\To V(N)$ is the qoutient map. The independence of the
expression \ref{bott} from the extensions $\talg{\alpha}$ and $\talg{s}$, as the properties of covariant derivation easily follow from the
definition by analogy to the case of the dual definition.

In general, if $\nabla$ is a covariant derivation on some vector bundle, its dual on the dual vector bundle is
defined by the equality
\begin{equation}\label{dualconns}
\scalar{\nabla_\alpha(\omega)}{s}=\alpha\scalar{\omega}{s}-\scalar{\omega}{\nabla_\alpha(s)}
\end{equation}
where $\omega$ and $s$ are sections of the bundles dual to each other. Keeping in mind this equality and follow the
definitions of the Bott connections on $\ortho$ and its dual $V(N)$, we obtain
$$
\lder{\zxa}{\talg{\omega}}(\talg{s})=\talg{s}\pars{\talg{\omega}(\zxa)}+
\zxa\pars{\talg{\omega}(\talg{s})}-\talg{s}\pars{\talg{\omega}(\zxa)}-
\talg{\omega}([\zxa,\talg{s}])=\zxa(\talg{\omega}(\talg{s}))-\talg{\omega}([\zxa,\talg{s}])
$$
from which follows the equality \ref{dualconns} for $\omega$ -- a section of $\ortho$, $s$ -- a section of $V(N)$
and $\alpha\in\restr$.

Let us consider an example which is rather simple but useful for the demonstration of various interesting
properties of singular Poisson structures.
\begin{exmp}\label{example}
Let $S$ be a 2-dimensional symplectic manifold, with Poisson bracket $\brck{\ }{\ }$. Let $\varphi$ be a smooth
function on $S$, such that $\varphi^{-1}(0)=\set{x_0}$. Conider a modified Poisson bracket
$\brck{\ }{\ }_1=\varphi\cdot\brck{\ }{\ }$ (as the manifold is 2-dimensional this bracket satisfies all the
properies required for a Poisson bracket). The symplectic foliation for the modified bracket consists of two
symplectic leaves: $\set{x_0}$ and $S\setminus\set{x_0}$. A fiber of the transversal foliation for the leaf
$S\setminus\set{x_0}$ is just $\set{0}$. Hence it is more interesting to consider the Bott connection for the
leaf consisting of just one point $\set{x_0}$. In this case the transversal space is $V(x_0)=T_{x_0}(S)$ and thus,
by definition, the covariant derivation, corresponding to the Bott connection, is an action of the elements of
$T_{x_0}^*(S)$ on $T_{x_0}(S)$:
\begin{equation}\label{bott1}
\nabla_\alpha(V)=[X^1_f,\talg{V}]|_{x_0}=[\varphi\cdot X_f,\talg{V}]|_{x_0}=V(\varphi)\cdot X_f(x_0)
\end{equation}
where $f\in\smooth{S}$ is such that $df(x_0)=\alpha$; $X^1_f$ is the Hamiltonian vector field for $f$
with respect to the modified bracket $\brck{\ }{\ }_1$ and $X_f$ is the Hamiltonian vector field for $f$
with respect to the original bracket $\brck{\ }{\ }$. Notice that as the original bracket is nondegenerate,
we have that $\alpha\neq0\tolfas X_\alpha\neq0$, therefore the equation $\nabla_\alpha(V)=0$, for $V$, when
$\alpha\neq0$ is equivalent to $V\in\ker(\varphi'(x_0))$. In other words, the ``flat'' sections of the vector
bundle $V(x_0)$ over $\set{x_0}$ are the elements of $\ker(\varphi'(x_0))$. Also it follows from the equality
\ref{bott1} that the Bott connection on the leaf $\set{x_0}$ is flat if and only if $\varphi'(x_0)=0$.
\end{exmp}
%%%%%%%%%%%%%%%%%%%%%%%%%%%%%%%%%%%%%%%%%%%%%%%%%%%%%%%%%%%%%%%
\section{Bott Connection for Transversal Mulivectors and Schouten-Nijenhuis Bracket}\label{bott-for-trans-multi}
For an arbitrary smooth manifold M we denote by $\fields{k}{M},\ k=0,\ldots,\infty$, the space of skew-symmetric
contravariant tensor (multivector) fields of degree $k$ on $M$ ($\fields{0}{M}=\smooth{M}$).
The Schouten-Nijenhuis bracket is a linear operation
$$
[\ ,\ ]:\fields{p}{M}\times\fields{q}{M}\To\fields{p+q-1}{M}
$$
with the following properties (see \cite{Nijen})
\begin{equation}\label{schout-props}
\begin{array}{ll}
(1) & [U,V]=(-1)^{|U|\cdot|V|}\cdot[V,U] \\
(2) & [U,V\wedge W]=[U,V]\wedge W+(-1)^{(|U|+1)\cdot|V|}\cdot V\wedge[U,W] \\
(3) & (-1)^{|U|\cdot(|W|-1)}\cdot[U,[V,W]]+(-1)^{|V|\cdot(|U|-1)}\cdot[V,[W,U]]+ \\
    & +(-1)^{|W|\cdot(|V|-1)}\cdot[W,[U,V]]=0
\end{array}
\end{equation}
For monomial type multivector fields the formula for the Schoute-Nijenhuis bracket is
\begin{equation}\label{schout}
\begin{array}{c}
[X_1\wedge\cdots\wedge X_m,Y_1\wedge\cdots\wedge Y_n]=\\
\\
(-1)^{m+1}\jami{i,j}{}(-1)^{i+j}[X_i,Y_j]\wedge X_1\wedge\cdots\wedge\hat{X_i}\wedge\cdots\wedge X_m\wedge Y_1\wedge\cdots\wedge\hat{Y_j}\wedge\cdots\wedge Y_n
\end{array}
\end{equation}
If we define a generalization of the ordinary Lie derivation for the multivector fields as
$$
L_X=i_X\circ d-(-1)^{|X|}d\circ i_X
$$
then there is the following relation between this operation and the Schouten-Nijenhuis bracket
\begin{equation}\label{lieschout}
[L_X,i_Y]=i_{[X,Y]}
\end{equation}
where $X$ and $Y$ are multivector fields and $i$ denotes the inner product of a multivector field and a differential
form.

For an arbitrary vector bundle $E$ and a covariant derivation $\nabla$ on it, the canonical extension of $\nabla$
to the Grassmann algebra bundle $\wedge(E)=\jami{k=0}{\infty}\wedge^k(E)$ is defined by the formula
\begin{equation}\label{covderext}
\nabla_X(S_1\wedge\cdots\wedge S_k)=\sum(-1)^{i+1}\nabla_X(S_i)\wedge S_1\wedge\cdots\wedge\hat{S_i}\wedge\cdots\wedge S_k
\end{equation}
where $S_1\wedge\cdots\wedge S_k$ is a section of $\wedge^k(E)$.

For a symplectic manifold $M$ and a symplectic leaf $N\subset M$, the $k$-th exterior degree of the vector bundle
$V(N)=\frac{T(M)|_N}{T(N)}$ is canonically isomorphic to $\frac{\wedge^kT(M)|_N}{T(N)\wedge(\wedge^{k-1}T(M))|_N}$.
Here $T(N)\wedge(\wedge^{k-1}T(M))|_N$ is the intersection of the ideal generated by $T(N)$ in the Grassmann algebra
bundle $\wedge T(M)|_N$, with $\wedge^kT(M)|_N$. According to the formulas \ref{schout} and
\ref{covderext}, the extension of the Bott covariant derivation on $V(N)$ to $\wedge^kV(N)$ can be defined as
\begin{equation}\label{BottSchouten}
\nabla_\alpha(U)=\rho([\zxa,\talg{U}]|_N)
\end{equation}
where $\alpha$ is a section of $\restr$; $U$ is a section of $\wedge^kV(N)$;
$$
\rho:\wedge^kT(M)|_N\To\wedge^kV(N)
$$
is the quotient map; $\talg{\alpha}$ is an extension of $\alpha$;
$\talg{U}$ is such a section of $\wedge^kT(M)$ (i.e., a multivector field on $M$) that $\rho(\talg{U}|_N)=U$
and $[\ ,\ ]$ is the \emph{Schouten-Nijenhuis Bracket}.
\begin{exmp}\label{example2}
Let $S$, $\brck{\ }{\ }_1=\varphi\brck{\ }{\ }$ and $x_0$ be the same objects as in the Example \ref{example}. We
have that $\wedge^2V(x_0)=\wedge^2T_{x_0}(S)$. For $\alpha\in T^*_{x_0}(S)$ and $U\wedge V\in T_{x_0}(S)$, we have
$$
\begin{array}{l}
\nabla_\alpha(U\wedge V)=\nabla_\alpha(U)\wedge V-\nabla_\alpha(V)\wedge U=\\
\\
U(\varphi)\cdot X_\alpha\wedge V-V(\varphi)\cdot X_\alpha\wedge U=X_\alpha\wedge(U(\varphi)\cdot V-V(\varphi)\cdot U)
\end{array}
$$
After this suppose that we need to solve the equation $\nabla_\alpha(U\wedge V)=0$ for $U\wedge V$. As
$\dim(S)=2$, the space $\wedge^2T_{x_0}(S)$ is 1-dimensional. Therefore we are free to take the vector $V$ as an
element of $\ker(\varphi'(x_0))$. Hence we obtain
$$
\nabla_\alpha(U\wedge V)=U(\varphi)\cdot X_\alpha\wedge V=0
$$
From this follows that if $\varphi'(x_0)=0$, then the extended Bott connection on $\wedge^2T_{x_0}(S)$ is flat,
otherwise, the equation $\nabla_\alpha(U\wedge V)=0$ has a nontrivial solution if and only if $X_\alpha\in\ker(\varphi'(x_0))$.
The both case are summarized as the following: the equation $\nabla_\alpha(U\wedge V)=0$ for has nontrivial solution if
and only if $X_\alpha\in\ker(\varphi'(x_0))$, and the set of solutions is $U\wedge\ker(\varphi'(x_0))$
(of course $\cong\wedge^2T_{x_0}(S)$).
\end{exmp}
%%%%%%%%%%%%%%%%%%%%%%%%%%%%%%%%%%%%%
\section{Generalized Functions (Distributions) on Poisson Manifold}
Traditionally, the generalized functions (distributions) on a finite-dimensional vector space are defined as
the elements of the topological dual to the space of test-objects which, itself, consists of the compactly
supported smooth functions. In this case, the space of test-objects can embedded into the space of distributions
by using of any volume form on the vector space. If we consider an arbitrary manifold, which is not necessarily
orientable, we should decide which set to use for test-objects. Besides the already mentioned space of functions
it can be the space of $n$-forms with compact support, where $n$ is the dimension of the manifold. In the case
of oriented manifold these spaces are isomorphic. In the case when the underlying manifold is not orientable
some constructions for these two realizations differ. For simplicity, we assume that $M$ is a compact, closed manifold
and first consider the case when the space of test-objects is $\smooth{M}$. We follow the standard notation and
denote the space of distributions on $M$ by $\distr{M}$. Usually we shall use the scalar product notation
$\scalar{\Phi}{f}$ for the value of $\Phi\in\distr{M}$ on $f\in\smooth{M}$.

The space \distr{M} is a module over the algebra \smooth{M}:
$$
\scalar{f\cdot\Phi}{g}=\scalar{\Phi}{fg},\quad f,g\in\smooth{M},\ \Phi\in\distr{M}
$$
It is also a Lie module over the Lie algebra of vector fields on $M$:
$$
\scalar{X(\Phi)}{f}=-\scalar{\Phi}{X(f)},\quad X\in\vecf{M}
$$
These two structures are correlated as
\begin{equation}\label{mult-der-correlation}
X(f\cdot\Phi)=X(f)\cdot\Phi+f\cdot X(\Phi)=(fX)(\Phi)
\end{equation}
\emph{Notice that in general, the equality $(fX)(\Phi)=f\cdot X(\Phi)$, is not true!}

If $M$ is oriented by a volume form $W$, define a scalar product on \smooth{M} by the standard formula
$$
\scalar{f}{g}=\int\limits_Mfg\cdot W,\quad f,g\in\smooth{M}
$$
We have an embedding
$$
\embm{W}:\smooth{M}\To\distr{M},\quad\emb{W}{f}=\scalar{f}{\ \cdot\ },\ \all f\in\smooth{M}
$$
\embm{W} is a homomorphism of \smooth{M}-modules:
$$
\emb{W}{fg}=f\cdot\emb{W}{g}\quad(\ =\emb{W}{f}\cdot g)
$$
therefore it is uniquely defined by its value on the constant function $1$. In general, \embm{W} is not a
homomorphism of $\vecf{M}$-Lie modules.
\begin{prop}\label{invariant-form-and-embedding}
For a given vector field $X\in\vecf{M}$, the following three conditions are equivalent
\begin{itemize}
\item[(1)]
	$\emb{W}{X(f)}=X(\emb{W}{f}),\quad\all f\in\smooth{M}$
\item[(2)]
	$X(\emb{W}{1})=0$
\item[(3)]
	the volume form $W$ is invariant under the 1-parameter flow of diffeomorphisms corresponding to $X$:
	$\lder{X}{W}=0$
\end{itemize}
\end{prop}
\begin{proof}
The implication $(1)\gamomd(2)$ is clear, and is obtained by taking $f\equiv1$ in (1). To show
$(2)\gamomd(1)$, recall that \embm{W} is a \smooth{M}-homomorphism, and thus $f\cdot\emb{W}{1}=\emb{W}{f}$. The
differentiation of the last equality by the vector field $X$, together with the relation
\ref{mult-der-correlation} gives
$$
X\pars{f\emb{W}{1}}=X(f)\cdot\emb{W}{1}=\emb{W}{X(f)}\gamomd X\pars{\emb{W}{f}}=\emb{W}{X(f)}
$$
As it is visible, $(1)\tolfas(2)$ is true for any \smooth{M}-homomorphism
$$
\embm{}:\smooth{M}\To\distr{M}
$$
$(2)\gamomd(3)$ can be obtained by direct calculation as follows
$$
\begin{array}{l}
0=\scalar{X\pars{\emb{W}{1}}}{f}=-\scalar{\emb{W}{1}}{X(f)}=-\int\limits_MX(f)\cdot W=\\
=-\int\limits_M\lder{X}{f\cdot W}+\int\limits_Mf\cdot\lder{X}{\omega}=-\int\limits_Md(f\cdot i_X(W))+\int\limits_Mf\cdot\lder{X}{W}=\\
=\int\limits_Mf\cdot\lder{X}{W},\ \all f\in\smooth{M}\gamomd\lder{X}{W}=0
\end{array}
$$
And $(3)\gamomd(2)$ can be obtained by following the above sequence of implications in the inverse order.
\end{proof}
Now assume that $M$ is a Poisson manifold. Define a binary operation
$$
\begin{array}{l}
\brck{\ }{\ }:\smooth{M}\times\distr{M}\To\distr{M},\quad\brck{f}{\Phi}=X_f(\Phi)\\
\\
f\in\smooth{M},\ \Phi\in\distr{M}
\end{array}
$$
which we regard as some kind of Poisson bracket of a smooth function and a distribution. This operation is
biderivative in the sense that it has the following properties
\begin{itemize}
\item[(1)]
	$\brck{f}{g\Phi}=\brck{f}{g}\Phi+g\brck{f}{\Phi}$
\item[(2)]
	$\brck{fg}{\Phi}=f\brck{g}{\Phi}+\brck{f}{\Phi}g$
\end{itemize}
Both of these are results of the relation \ref{mult-der-correlation}: (1) is obtained directly from
$X(f\Phi)=X(f)\Phi+fX(\Phi)$ and (2) from $(fX)(\Phi)=X(f)\Phi+fX(\Phi)$ as follows
$$
\begin{array}{c}
\brck{fg}{\Phi}=(fX_g+gX_f)(\Phi)=X_g(f)\cdot\Phi+f\cdot X_g(\Phi)+X_f(g)\cdot\Phi+g\cdot X_f(\Phi)=\\
=(\underbrace{\brck{g}{f}+\brck{f}{g}}_0)\cdot\Phi+f\brck{g}{\Phi}+g\brck{f}{\Phi}=f\brck{g}{\Phi}+g\brck{f}{\Phi}
\end{array}
$$
Under these conditions we state that \distr{M} is a \emph{Poisson module} over the Poisson algebra \smooth{M}. The
formal definition of a Poisson module is the following
\begin{defn}[Poisson Module]
Let $A$ be an associative Poisson algebra. A space $V$ is said to be a Poisson module over $A$ if $V$ is an
$(A,A)$-bimodule, when $A$ is considered as an associative algebra; $V$ is a Lie module over $A$, when $A$ is
considered as a Lie algebra (and we denote the action of $a\in A$ on $v\in V$ by $\brck{a}{v}$); for any $a,b\in A$
and $v\in V$, the following equalities are true
$$
\begin{array}{l}
\brck{ab}{v}=a\brck{b}{v}+\brck{a}{v}b\\
\brck{a}{bv}=b\brck{a}{v}+\brck{a}{b}v\textrm{ and }\brck{a}{vb}=\brck{a}{v}b+v\brck{a}{b}
\end{array}
$$
\end{defn}
It is clear that for a Poisson manifold $M$, the space \smooth{M} is a Poisson module over itself. Thus it is
natural to look for such embeddings of \smooth{M} into \distr{M} that are homomorphisms of \smooth{M}-Poisson
modules.
%%%%%%%%%%%%%%%%%%%%%%%%%%%%%%%%%%%%%%%%%%%%%%%
\section{Hamiltonian-Invariant Volume Forms and the Modular Class of a Poisson Manifold: a Brief Review}
As it is well-known (see \cite{Koszul}, \cite{Vaisman}, \cite{Lichnerowich}), the bivector field $P$ corresponding
to the Poisson structure on the manifold $M$ defines an operator
$$
\begin{array}{l}
\sigma:\fields{k}{M}\To\fields{k+1}{M},\quad\sigma(W)=[P,W],\ k=1,\ldots,\infty\\
\\
\sigma(f)=X_f,\textrm{ for }f\in\fields{0}{M}=\smooth{M}
\end{array}
$$
This operator, due to the properties \ref{schout-props} of the Schouten-Nijenhuis bracket, is coboundary and
antiderivative. Therefore it defines a cohomology algebra, which is known as the \emph{Poisson cohomology} of the
Poisson algebra \smooth{M}. As it follows directly from the definition, in the dimension 1 the space $\operatorname{Im}(\sigma)$
is the space of Hamiltonian vector fields and $\ker(\sigma)$ is the space of such vector fields, the 1-parameter
flows of which preserve the Poisson bracket.

\newcommand{\mcl}[1]{\mu_{_{#1}}}
The \emph{modular class} of the Poisson manifold $M$ is known as the obstruction to the existence of a volume form on
$M$ which is invariant for the Hamiltonian flows (see \cite{Weinstein}, \cite{Ginzburg}). Let us give a brief
review of this relation. First of all notice that the Poisson manifold must be orientable. Assume that $W$ is any
volume form on $M$. Consider an operator $\mcl{W}:\smooth{M}\To\smooth{M}$ defined by the equality
$$
\mcl{W}(f)\cdot W=\lder{X_f}{W},\quad\all f\in\smooth{M}
$$
This operator is a derivation: for $f,g\in\smooth{M}$ we have the following
$$
\begin{array}{c}
X_{fg}=fX_g+gX_f\gamomd\mcl{W}(fg)\cdot W=\lder{fg}{W}=\\
\\
=df\wedge(i_{X_g}W)+dg\wedge(i_{X_f}W)+
f\underbrace{\lder{X_g}{W}}_{\mcl{W}(g)\cdot W}+g\underbrace{\lder{X_f}{W}}_{\mcl{W}(f)\cdot W}=\\
\\
=-i_{X_g}(\underbrace{df\wedge W}_0)+\brck{g}{f}\cdot W-i_{X_f}(\underbrace{dg\wedge W}_0)+\brck{f}{g}\cdot W+\\
\\
+(f\mcl{W}(g)+g\mcl{W}(f))\cdot W=(f\mcl{W}(g)+g\mcl{W}(f))\cdot W
\end{array}
$$
Thus $\mcl{W}$ can be considered as a vector field on $M$. It is known as the \emph{modular vector field} corresponding
to the volume form $W$. If $W_1$ is another volume form, then $W_1=\varphi\cdot W,\ \varphi\in\smooth{M},\ \varphi^{-1}(0)=\emptyset$.
The relation between $\mcl{W}$ and $\mcl{W_1}$ is obtained as follows
$$
\begin{array}{c}
\mcl{W_1}(f)\cdot W_1=\lder{X_f}{W_1}=\brck{f}{\varphi}\cdot W+\varphi\cdot\mcl{W}(f)\cdot W=\\
\\
=(\mcl{W}(f)-\frac{X_\varphi(f)}{\varphi})\cdot W_1=(\mcl{W}(f)-X_{\ln(\varphi)})\cdot W_1\gamomd\\
\\
\gamomd\mcl{W}-\mcl{W_1}=\sigma(\ln(\varphi))
\end{array}
$$
For any $f,g\in\smooth{M}$ we have
$$
\begin{array}{c}
\mcl{W}(\brck{f}{g})\cdot W=\lder{X_{\brck{f}{g}}}{W}=[L_{X_f},L_{X_g}](W)=\\
\\
=(X_f(\mcl{W}(g))-X_g(\mcl{W}(f)))\cdot W=(\brck{\mcl{W}(f)}{g}+\brck{f}{\mcl{W}(g)})\cdot W
\end{array}
$$
which is the infinitesimal version of the fact that the 1-parameter flow of $\mcl{W}$ preserves the Poisson bracket:
$\sigma(\mcl{W})=0$. The latter means that any modular vector field is a cocycle in the complex $(\fields{}{M},\sigma)$,
and the relation between $\mcl{W}$ and $\mcl{W_1}$ implies that all modular vector fields define one and the same
cohomology class in the Poisson cohomology algebra. Let us denote this class by $\mu$. This cohomology class is
know as the \emph{modular class} of the Poisson structure. If there exists a Hamiltonian invariant volume form
$W_1=\varphi\cdot W$, then $\mcl{W_1}=0$ and from the relation between $\mcl{W}$ and $\mcl{W_1}$ we obtain that
$\mcl{W}=\sigma(\ln(\varphi))$, which implies that the modular class is trivial. Conversely, if the modular class
is trivial, then we have
$$
\begin{array}{c}
\mcl{W}=\sigma(\psi)=X_\psi,\textrm{ for some }\psi\in\smooth{M}\gamomd\lder{X_f}{W}=\brck{\psi}{f}\cdot W\gamomd\\
\\
\gamomd\lder{X_f}{W}+\brck{f}{\psi}\cdot W=0\gamomd\brck{f}{\psi}\cdot\exp(\psi)\cdot W+\exp(\psi)\cdot\lder{X_f}{W}=\\
\\
=\lder{X_f}{\exp(\psi)\cdot W}=0,\quad\all f\in\smooth{M}
\end{array}
$$
That is, the volume form $\exp(\psi)\cdot W$ is invariant for the Hamiltonian flows on $M$.
\begin{exmp}\label{example-nontrivial-modular}
Consider the 2-dimensional Poisson manifold $(S,\brck{\ }{\ }_1=\varphi\brck{\ }{\ })$ that we introduced in the
Example \ref{example}. Notice that for any symplectic manifold its modular class is trivial, because at least one
(and only one up to a constant multiplier) Hamiltonian invariant volume form is $\omega^n$, where $\omega$ is the
symplectic form and $2n$ is the dimension of the manifold. Hence the modular class of the original Poisson bracket
on $S$ is trivial. Let us calculate the modular class with respect to the deformed Poisson bracket
$\brck{\ }{\ }_1=\varphi\cdot\brck{\ }{\ }$, which is singular in the point $x_0=\varphi^{-1}(0)$. As a volume form
we can take the original symplectic form $\omega$. Then we have
$$
\begin{array}{c}
\lder{X^1_f}{\omega}=\lder{\varphi X_f}{\omega}=d\varphi\wedge(i_{X_f}\omega)+\varphi\cdot\underbrace{\lder{X_f}{\omega}}_0=\\
=-i_{X_f}(\underbrace{d\varphi\wedge\omega})+\brck{f}{\varphi}\cdot\omega=X_\varphi(f)\cdot\omega
\end{array}
$$
Thus we obtain that the modular vector field for $\omega$, with respect to the modified Poisson structure is the
Hamiltonian vector field of $\varphi$, but with respect to the original Poisson bracket. To be cohomological to
$0$, this vector field must be Hamiltonian with respect to the modified Poisson bracket too:
$$
\begin{array}{c}
X_\varphi=X^1_f\tolfas X_\varphi=\varphi\cdot X_f\gamomd\frac{X_\varphi}{\varphi}=X_f,\textrm{ on }S\setminus\set{x_0}\tolfas\\
\tolfas X_{\ln(\varphi)}=X_f\textrm{ on }S\setminus\set{x_0}\tolfas f-\ln(\varphi)=\textbf{const}\textrm{ on }S\setminus\set{x_0}
\end{array}
$$
As the point $x_0$ is singular for $\ln(\varphi)$, there is no such $f\in\smooth{S}$ that $f-\ln(\varphi)=\textbf{const}$.
Therefore the modular class for $(S,\brck{\ }{\ }_1)$ is not trivial. We can conclude that there is no such volume
form on $S$ which is invariant for the Hamiltonian flows with respect to the bracket $\brck{\ }{\ }_1=\varphi\cdot\brck{\ }{\ }$.
\end{exmp}
%%%%%%%%%%%%%%%%%%%%%%%%%%%%%%%%%%%%%%%%%
\section{The Generalized Center of a Poisson Structure}
As it follows from the Proposition \ref{invariant-form-and-embedding}, if $M$ is a Poisson manifold and $W$ is
a volume form on it, the embedding $\embm{W}:\smooth{M}\To\distr{M}$ is a homomorphism of Poisson modules if and
only if $W$ is an invariant form for all Hamiltonian flows. But for degenerate Poisson structures the existence of
such volume form, in general, is problematic (see Example \ref{example-nontrivial-modular}). Some kind of
alternatives to invariant volume forms can be the generalized central elements (Casimirs) for the Poisson structure.
\begin{defn}
We call the \emph{generalized center} of the Poisson algebra \smooth{M} the subspace of \distr{M} defined as
$$
\set{\Phi\in\distr{M}\ |\ \brck{f}{\Phi}=0,\ \all f\in\smooth{M}}\equiv\distro{M}
$$
\end{defn}
It follows from the definition that a \smooth{M}-homomorphism
$$
\mathcal{E}:\smooth{M}\To\distr{M}
$$
is a homomorphism of the Poisson modules if and only if $\mathcal{E}(1)\in\genc{M}$.

Let us consider \genc{M} from the homological viewpoint. There is a well-defined boundary operator of degree -1
on the exterior algebra of differential forms on the Poisson manifold $M$ (see \cite{Brylinski}, \cite{Vaisman})
$$
\delta:\forms{\bullet}{M}\To\forms{\bullet-1}{M},\ \delta=i_P\circ d-d\circ i_P
$$
where $P$ is the bivector field corresponding to the Poisson structure. The ''expanded`` definition of $\delta$
is the following
$$
\begin{array}{c}
\delta(\phi_0d\phi_1\wedge\cdots\wedge\phi_k)=\jami{i=1}{k}(-1)^{i+1}\brck{\phi_0}{\phi_i}d\phi_1\wedge\cdots\wedge\hat{d\phi_i}\wedge\cdots\wedge d\phi_k+\\
+\jami{i<j}{}(-1)^{i+j}\phi_0d\brck{\phi_i}{\phi_j}\wedge d\phi_1\wedge\cdots\wedge\hat{d\phi_i}\wedge\cdots\wedge\hat{d\phi_j}\wedge\cdots\wedge d\phi_k
\end{array}
$$
The complex $(\allforms{M},\delta)$ is known as the \emph{canonical complex of the Poisson manifold} $M$, and its
homology space $H_*(M,\delta)$ --- the \emph{canonical homology of the Poisson structure on} $M$. Assume that $M$ is a symplectic manifold, or
equivalently, the bivector field $P$ is nondegenerate. $P$ can be considered as an antisymmetric \smooth{M}-linear
mapping
$$
P:\forms{1}{M}\times\forms{1}{M}\To\smooth{M}
$$
which can be extended to higher degrees
$$
P^k:\forms{k}{M}\times\forms{k}{M}\To\smooth{M},\quad k=1,\ldots,\infty
$$
The symplectic analogue to the star operator for a Riemann manifold is defined as (see \cite{Brylinski})
$$
\begin{array}{l}
*:\forms{k}{M}\To\forms{2n-k}{M},\quad\omega_1\wedge(*\omega_2)=P^k(\omega_1\wedge\omega_2)\cdot\frac{\omega^n}{n!}\\
\\
\textrm{for }\omega_1,\omega_2\in\forms{k}{M},\ k=1,\ldots,\infty
\end{array}
$$
where $2n=\dim(M)$ and $\omega$ is the symplectic form on $M$, corresponding to $P$. One of the properties of this
operator is the following relation with $\delta$
$$
\delta=(-1)^{k+1}*\circ d\circ *:\forms{k}{M}\To\forms{k-1}{M}
$$
This property implies that $*$ induces an isomorphism of the canonical homologies $H_k(M,\delta)$ and the De Rham
cohomologies $H^{2n-k}(M)$, for $k=0,\ldots,2n$. It follows from the definition of $\delta$ that
$$
\delta\pars{\forms{1}{M}}=\brck{\smooth{M}}{\smooth{M}}=\set{\jami{i}{}\brck{f_i}{g_i}\ |\ f_i,g_i\in\smooth{M}}
$$
Therefore we have
$$
H_0(M,\delta)=\frac{\smooth{M}}{\brck{\smooth{M}}{\smooth{M}}}
$$
The symplectic manifold $M$ is oriented by the volume form $\omega^n$. If $M$ is compact, then we have that
$H^{2n}(M)\cong\Real$, after which, because of the isomorphism in homologies induced by $*$,we obtain
\begin{equation}\label{commutator-quotient}
\frac{\smooth{M}}{\brck{\smooth{M}}{\smooth{M}}}=H_0(M,\delta)\cong\Real
\end{equation}
Now return back to the generalized Casimirs and notice that if $\Phi\in\genc{M}$, then by definition
$$
\scalar{\brck{f}{\Phi}}{g}=-\scalar{\Phi}{\brck{f}{g}}=0,\quad\all f,g\in\smooth{M}
$$
Hence we can state that \genc{M} is the dual linear space of $H_0(M,\delta)$
$$
\genc{M}=H_0(M,\delta)^*
$$
After which, from the isomorphism \ref{commutator-quotient} follows that in the case of symplectic manifold,
the space \genc{M} is 1-imensional and is spanned by the element \emb{{\omega^n}}{1}, where \embm{{\omega^n}} is
the embedding
$$
\embm{{\omega^n}}:\smooth{M}\To\distr{M},\quad\emb{{\omega^n}}{f}=\int\limits_Mf\cdot\omega^n
$$
This also implies that in the symplectic case, there is one and only one (up to a constant multiplier) homomorphism
of \smooth{M}-Poisson modules \smooth{M} and \distr{M}.

Now consider an example of a degenerate Poisson structure
\begin{exmp}\label{dirac-genc}
For the Poisson manifold $(S,\brck{\ }{\ }_1=\varphi\cdot\brck{\ }{\ })$ (see Example \ref{example}),
one element of the generalized center is the Dirac functional $\delta_{x_0}$. For this distribution we have
$$
\begin{array}{c}
\scalar{\brck{f}{\delta_{x_0}}_1}{g}=-\scalar{\delta_{x_0}}{\brck{f}{g}_1}=\\
\\
=-\scalar{\delta_{x_0}}{\varphi\brck{f}{g}}=-\varphi(x_0)\brck{f}{g}(x_0)=0,\quad\all f,g\in\smooth{S}
\end{array}
$$
As the space \genc{S} is not trivial, we have at least one homomorphism of Poisson modules
$$
\embm{\delta_{x_0}}:\smooth{S}\To\distr{S},\quad\emb{\delta_{x_0}}{f}=f\cdot\delta_{x_0}
$$
but it is clear that this is not an embedding:
$$
\ker(\embm{\delta_{x_0}})=I_{x_0}=\set{f\in\smooth{S}\ |\ f(x_0)=0}
$$
We have that for any Poisson manifold the ''ordinary`` center of the Poisson algebra of smooth functions is
always at least 1-dimensional and contains the set of all constant functions. Also we have that if this center is
more then 1-dimentional, the Poisson bracket is clearly degenerate. But as it follows from this example, the
converse is not true: because the bracket $\brck{\ }{\ }_1$ is degenerate in one point, the center of the
Poisson algebra \smooth{S} consists just the constants. We expect that in such cases the dimension of the space
of generalized Casimirs is $>1$. Though in the case of this example it is not so clear, because we found just
1-dimensional subspace of \genc{S}. Further we shall describe some other ways of construction of generalized
Casimir elements.
\end{exmp}
%%%%%%%%%%%%%%%%%%%%%%%%%%%%%%%%%%%%%%%%%
\section{Generalized Casimirs Supported by Symplectic Leaves}
For any smooth map $f:M_1\To M_2$, where $M_1$ and $M_2$ are smooth manifolds, we have the dual
map $f^*:\smooth{M_2}\To\smooth{M_1}$. The latter, itself, induces a map
$f_*:\distr{M_1}\To\distr{M_2}$: $\scalar{f_*(\Phi)}{\psi}=\scalar{\Phi}{f^*(\psi)}$. Assume that $M_1$ and
$M_2$ are Poisson manifolds and $f$ is a Poisson map, i.e.,
$f^*(\brck{\phi}{\psi}_2)=\brck{f^*(\phi)}{f^*(\psi)}_1$, where $\brck{\ }{\ }_1$
($\brck{\ }{\ }_2$) denotes the Poisson bracket on $M_1$ ($M_2$). Under these conditions, the
mapping $f_*:\distr{M_1}\To\distr{M_2}$ is a homomorphism of Poisson modules, in the sense that
$$
\brck{\psi}{f_*(\Phi)}=f_*(\brck{f^*(\psi)}{\Phi}),\quad\all\psi\in\smooth{M_2}\textrm{ and }\all\Phi\in D(M_1)
$$
As a result we obtain the following
\begin{lem}
If $f:M_1\To M_2$ is a Poisson map for the Poisson manifolds $M_1$ and $M_2$ then
$f_*:\distr{M_1}\To\distr{M_2}$ maps the generalized center of $\smooth{M_1}$ into the generalized center
of $\smooth{M_2}$.
\end{lem}
This construction can be applied to the case when $M_2\equiv M$ is a Poisson manifold,
$M_1\equiv N\subset M$ is a Poisson submanifold of $M$ and $f:N\To M$ is the embedding map.
In particular, consider the case when the Poisson bracket on $M$ is degenerate and $N\subset M$ is
a symplectic leaf. We assume that $N$ is closed and compact submanifold. As $N$ is a compact
symplectic manifold, the generalized center of $\smooth{N}$ is one-dimensional and is generated by
the functional $\emb{\omega_N}{1}$, where $\omega_N$ is the symplectic form on $N$ induced by the Poisson
structure on $N$. According to the above lemma, the mapping $J_*:\distr{N}\To\distr{M}$,
where $J:N\hookrightarrow M$ is the
embedding map, maps the functional $\emb{\omega_N}{1}$ into the generalized center of the Poisson
algebra $\smooth{M}$. Let us denote the image of $\emb{\omega_N}{1}$ into $\genc{M}$ by $\delta_N$
(this is some generalization of the Dirac's generalized function, which justifies this notation).
The explicit expression for $\delta_N$ is the following
\begin{equation}\label{leaf-delta-formula}
\delta_N(\phi)=\int\limits_N\phi|_{_N}\cdot\omega_{_N}^k,\qquad\all\phi\in\smooth{M}
\end{equation}
where $k=\frac{1}{2}\dim(N)$.
\begin{prop}\label{set_of_leaves}
For any finite set of compact closed symplectic leaves $N_1$, $\ldots$, $N_m$, the corresponding
functionals $\delta_{N_i},\,i=1,\ldots,m$ are linearly independent.
\end{prop}
\begin{proof}
Consider a set of points $x_i\in N_i,\,i=1,\ldots,m$. As the leaves $N_1$, $\ldots$, $N_m$ are
compact, for any $i$ we can select an open neighborhood of the point $x_i$: $U_i\subset M$, so
that $U_i\cap U_j=U_i\cap N_j=\emptyset$ when $i\neq j$. For each $U_i$ consider a function
$f_i\in\smooth{M}$ such that: $f_i\geq 0$, $f_i(x_i)>0$ and $f_i(x)=0$ when $x\notin U_i$. From
the formula \ref{leaf-delta-formula} easily follows that $\delta_{N_i}(f_j)=0$ when $i\neq j$ and
$\delta_{N_i}(f_i)\equiv a_i>0$. Assume that $\jami{i=1}{m}\alpha_i\cdot\delta_{N_i}=0$ for some
$\alpha_i\in\Real,\,i=1,\ldots,m$. The equalities
$$
(\jami{i=1}{m}\alpha_i\cdot\delta_{N_i})(f_j)=\alpha_j\cdot a_j=0,\quad j=1,\ldots,m
$$
imply that $\alpha_j=0$, $j=1,\ldots,m$.
\end{proof}
\begin{exmp}
In the case of the Poisson manifold $(S,\brck{\ }{\ }_1=\varphi\cdot\brck{\ }{\ })$ from the Example \ref{example},
for the symplectic leaf $\set{x_0}=\varphi^{-1}(0)$, the formula \ref{leaf-delta-formula} gives the Dirac functional
$\delta_{x_0}$ as it was already described in the example at the end of the previous section. Another symplectic
leaf is $S\setminus\set{x_0}$, which is not compact. The restricted Poisson structure on this leaf is nondegenerate
and the corresponding symplectic form is $\frac{\omega}{\varphi}$, where $\omega$ is the symplectic form on $S$,
corresponding to the original bracket $\brck{\ }{\ }$. Clearly this form is singular when $x\to x_0$ and therefore
when $N=S\setminus\set{x_0}$ we cannot construct an element of \genc{S} by the same method as $\delta_N$, for
a compact $N$. Though we can construct construct some other elements of \genc{S} as follows. Consider a non-zero
tangent vector $v\in T_{x_0}(S)$, such that $\varphi'_{x_0}(v)=0$. The distribution
$$
\delta'_{x_0}(v):\ \scalar{\delta'_{x_0}(v)}{f}=-f'_{x_0}(v),\ \all f\in\smooth{S}
$$
which is the derivative of the Dirac functional $\delta_{x_0}$ by $v$, is also an element of \genc{S}:
$$
\begin{array}{c}
\scalar{\brck{g}{\delta'_{x_0}(v)}_1}{f}=-\scalar{\delta'_{x_0}(v)}{\brck{g}{f}_1}=\scalar{\delta'_{x_0}(v)}{\varphi\brck{f}{g}}=\\
\\
=\underbrace{\varphi'_{x_0}(v)}_0\brck{f}{g}(x_0)+\underbrace{\varphi(x_0)}_0\brck{f}{g}'_{x_0}(v)=0,\quad\all f,g\in\smooth{M}
\end{array}
$$
It is clear that the distributions $\delta_{x_0}$ and $\delta'_{x_0}(v)$ are linearly independent and therefore
$\dim(\genc{S})\geq2$. If the point $x_0$ is singular for the function $\varphi$, i.e., $\varphi'_{x_0}\equiv0$, then
$\dim(\genc{S})\geq3$.
\end{exmp}
Let us generalize the construction described in the above example for any compact symplectic leaf of the Poisson
manifold $M$. Denote by $W(N)$ the space of such linear maps $X:\smooth{M}\To\smooth{N}$ that
$$
X(fg)=X(f)\rho(g)+\rho(f)X(g),\quad\all f,g\in\smooth{M}
$$
where $\rho:\smooth{M}\To\smooth{N}$ is the restriction map.\\
\textbf{Remark:} in fact, such $X$ is a section of the bundle $T(M)|_N$.\\
Any such $X$ defines a mapping (the dual of $X$) $\distr{N}\To\distr{M}$, which we denote also by $X$:
$$
\scalar{X(\Phi)}{f}=-\scalar{\Phi}{X(f)},\quad\Phi\in\distr{N},\ f\in\smooth{M}
$$
\begin{lem}
for $X\in W(N)$ the following two comditions are equivalent
$$
\begin{array}{ll}
(1) & X\pars{\brck{f}{g}}=\brck{X(f)}{\rho(g)}+\brck{\rho(f)}{X(g)},\ \all f,g\in\smooth{M}\\
&\\
(2) & X\pars{\brck{\rho(f)}{\Phi}}=\rho^*\pars{\brck{X(f)}{\Phi}}+\brck{f}{X(\Phi)},\ \all f\in\smooth{M},\ \all\Phi\in\distr{N}
\end{array}
$$
where $\rho^*:\distr{N}\To\distr{M}$ is the dual of the restriction map $\rho$.
\end{lem}
\begin{proof}
The second equality is for distributions. Hence if we rewrite the both sides of it in the terms of their values
on a function $g\in\smooth{M}$, we obtain
$$
\scalar{\Phi}{X\brck{f}{g}}=\scalar{\Phi}{\brck{X(f)}{\rho(g)}}+\scalar{\Phi}{\brck{\rho(f)}{X(g)}},\quad\all\Phi\in\distr{N}
$$
which implies (1) for $X,\ f$ and $g$. From the latter equality clearly follows $(1)\gamomd(2)$.
\end{proof}
Let us denote by $\talg{W}(N)$ the subspace of such elements $X\in W(N)$ which satisfy the identity (1). It follows
from the above lemma that for the elements of $\talg{W}(N)$, the corresponding duals $X:\distr{N}\To\distr{M}$ maps
the generalized center \genc{N} to the generalized center \genc{M}. If $X=X_\psi,\ \psi\in\smooth{N}$, is a Hamiltonian vector field on $N$,
then the identity (1) follows from the Jacoby identity for the Poisson bracket.
Therefore such $X_\psi$ is always an element of $\talg{W}(N)$. But it follows from the definition of \genc{N} that
$X_\psi\pars{\genc{N}}=\set{0}$ (because $X_\psi\pars{\genc{N}}=-\brck{\psi}{\Phi},\ \Phi\in\distr{N}$).
Therefore we can consider the quotient of $\talg{W}(N)$ by the linear subspace generated by the Hamiltonian vector
fields $X_\psi,\ \psi\in\smooth{N}$. Let us denote this quotient by $W_0(N)$.
\begin{prop}
For $X\in W(N)$, the following two statements are equivalent
\begin{itemize}
\item[(1)]
	$X(\delta_N)=0$
\item[(2)]
	$X$ is a tangent vector field on $N$, and its 1-parameter flow of diffeomorphisms preserves
	the volume form $\omega^k_N,\ k=\frac{1}{2}\dim(N)$
\end{itemize}
\end{prop}
\begin{proof}
To prove $(1)\gamomd(X\textrm{ is a tangent vector field on }N)$ it is sufficient to show $X(I_N)=\set{0}$, where
$I_N\subset\smooth{M}$ is the ideal of functions vanishing on $N$. For $f\in I_N$ we have
$$
\begin{array}{c}
\scalar{X(\delta_N)}{fg}=-\int\limits_N(X(f)\rho(g)+\underbrace{\rho(f)}_0X(g))\cdot\omega_N^k=\\
\\
=-\int\limits_NX(f)\rho(g)\cdot\omega_N^k=0,\quad\all g\in\smooth{M}\gamomd X(f)=0
\end{array}
$$
Now let $X$ be a tangent vector field on $N$ and $X(\delta_N)=0$. Then we have
$$
\begin{array}{c}
di_X\pars{\rho(f)\cdot\omega_N^k}=\lder{X}{\rho(f)\cdot\omega_N^k}=X(\rho(f))\cdot\omega_N^k+\rho(f)\cdot\lder{X}{\omega_N^k}\gamomd\\
\\
\gamomd0=\scalar{X(\delta_N)}{f}=\int\limits_N\rho(f)\cdot\lder{X}{\omega_N^k},\ \all f\in\smooth{M}\gamomd\lder{X}{\omega_N^k}=0
\end{array}
$$
Hence the implication $(1)\gamomd(2)$ is proved. The converse easily follows from the last serie of equalities and
the implications taken in the reverse order.
\end{proof}
To summarize, we can say that besides the generalized Casimirs of the type $\delta_N$, for compact symplectic
leaves $N$, we have its ''derivatives`` --- $X(\delta_N)$, by the elemements $X\in W_0(N)$.

There is a canonical mapping $W(N)\To\Gamma\pars{V(N)}$, where $V(N)$ is the transversal vector bundle over $N$ used
in the definition of the Bott connection: $V(N)=\frac{T(M)|_N}{T(N)}$, and $\Gamma$ denotes the space of sections.
Clearly, the subspace of Hamiltonian fields on $N$ is mapped to $\set{0}$. Thus we have a well-defined mapping
$\pi:W_0(N)\To V(N)$. The condition $U\brck{f}{g}=\brck{U(f)}{\rho(g)}+\brck{\rho(f)}{U(g)}$ for $U\in W(N)$ and
$f,g\in\smooth{M}$, rewritten as
$$
\pars{\pars{U\circ X_f-X_f\circ U}g}|_N=X_{U(f)}\pars{\rho(g)}\quad\pars{\ \tolfas-[X_f,U]|_N=X_{U(f)}}
$$
implies, by definition of the Bott connection on $V(N)$, the equality
$$
\nabla_{df}\pars{\pi(U)}=0,\quad\all f\in\smooth{N}
$$
Therefore the image of the mapping $\pi:W_0(N)\To V(N)$ is the set of flat sections of $V(N)$, with respect to
the Bott connection on it.

The support set of a generalized Casimir is foliated by some subset of the set of symplectic leaves of the
Poisson manifold. The following will shed a light on the meaning of this statement.

For any open subset $U\subset M$ let us denote by $\smooth{M}_U$ the subspace of $\smooth{M}$
consisting of the functions vanishing on $M\setminus U$. The restriction of a distribution $u\in\distr{M}$
to $U$ is defined as its restriction to the subspace $\smooth{M}_U$.
Let us recall the definition of the support set of a generalized function.
\begin{defn}[The support set of a distribution, see \cite{Hormander}]
For $u\in\distr{M}$, the support of $u$, denoted by $Supp(u)$, is the set of points in $M$ having no
open neighborhood to which the restriction of $u$ is 0.
\end{defn}
\begin{defn}[Singular support, see \cite{Hormander}]
For $u\in\distr{M}$, the singular support of $u$ is the set of points in $M$ having no open
neighborhood to wich the restriction of $u$ is a $C^\infty$ function.
\end{defn}
If $f:M_1\To M_2$ is diffeomorphism of the manifolds $M_1$ and $M_2$, then for any $\Phi\in\distr{M_1}$
we have that $Supp(f_*(\Phi))=f(Supp(\Phi))$. Therefore, if we have an action of a Lie group $G$
on the manifold $M$, and a generalized function $\Phi\in\distr{M}$ is invariant under this action,
then its support set is also an invariant set for this action. Moreover, if $S\subset M$
is such that $S$ is invariant set for the action of $G$, the action on $S$ is transitive and
$Supp(\Phi)\subset S$ for some invariant $\Phi\in\distr{M}$, then $Supp(\Phi)=S$. In this case, for
any point $x\in Supp(\Phi)$, the set $Supp(\Phi)$ is the orbit of $x$ under the action of the
group $G$.
\begin{lem}
Let $\set{\varphi_t}_{t\in\Real}$ be the 1-parameter group of diffeomorphisms of $M$ corresponding
to a Hamiltonian vector field $Ham(f)$ for some $f\in\smooth{M}$. If $\Phi$ is an element of \genc{M},
then $\Phi$ is invariant under the action of the group $\set{\varphi_t}_{t\in\Real}$.
\end{lem}
\begin{proof}
By definition we have $\scalar{{\varphi_*}_t(\Phi)}{g}=\scalar{\Phi}{g\circ\varphi^{-1}_t}$,
$g\in\smooth{M}$. This together with the definition of the generalized center gives:
$\frac{d}{dt}\scalar{\Phi}{g\circ\varphi^{-1}_t}=\scalar{\Phi}{-\brck{f}{g}}=\scalar{\brck{f}{\Phi}}{g}=0$,
$\all g\in\smooth{M}$.
\end{proof}
\begin{cor}
If $\Phi$ is an element of \genc{M} and $N\cap Supp(\Phi)\neq\emptyset$, for some compact and closed symplectic leaf
$N\subset M$, then $N\subset Supp(\Phi)$.
\end{cor}
%%%%%%%%%%%%%%%%%%%%%%%%%%%%%%%%%%%%%%%%%%%%%%
\section{Another Realization of the Space of Distributions on a Poisson Manifold}
Another realization of the space of distributions is obtained when as the space of test-objects is taken the space
\forms{n}{M} --- the space of n-differential forms on $M$ where $n=\dim(M)$. Hence, in this case, the space of
distributions on $M$ is the topological dual to \forms{n}{M}. We denote this space by \distrf{M}. As \distr{M},
this space also is a module over \smooth{M}:
$$
\scalar{f\cdot\Phi}{\alpha}=\scalar{\Phi}{f\alpha},\quad\Phi\in\distrf{M},\ f\in\smooth{M}
$$
And is a Lie module over $\vecf{M}$ under the following action of a vector field
$$
\scalar{X(\Phi)}{\alpha}=-\scalar{\Phi}{\lder{X}{\alpha}},\quad X\in\vecf{M},\ \alpha\in\forms{n}{M}
$$
Unlike the space \distr{M}, in this case we have a canonical embedding
$$
\mathcal{E}:\smooth{M}\To\distrf{M}
$$
defined as
$$
\scalar{\mathcal{E}(f)}{\alpha}=\int\limits_Mf\cdot\alpha,\quad f\in\smooth{M},\ \alpha\in\forms{n}{M}
$$
Clearly $\mathcal{E}$ is a homomorphism of \smooth{M}-modules, and again unlike the case of \distr{M}, it is a
homomorphism of $\vecf{M}$-Lie modules
$$
\begin{array}{c}
\scalar{\mathcal{E}\pars{X(f)}}{\alpha}=\int\limits_MX(f)\cdot\alpha=\int\limits_M\lder{X}{f\cdot\alpha}-\int\limits_Mf\cdot\lder{X}{\alpha}=\\
\\
=\int\limits_Md\pars{i_X\pars{f\cdot\alpha}}-\int\limits_Mf\cdot\lder{X}{\alpha}=-\int\limits_Mf\cdot\lder{X}{\alpha}=\scalar{X\pars{\mathcal{E}(f)}}{\alpha}\\
\\
\all X\in\vecf{M},\ \all f\in\smooth{M}\textrm{ and }\all\alpha\in\forms{n}{M}
\end{array}
$$
Therefore, if $M$ is a Poisson manifold (and hence \distrf{M} is a Poisson module over \smooth{M}), the embedding
$\mathcal{E}$ is a homomorphism of Poisson modules. This implies that the center of the Poisson algebra \smooth{M}
is mapped to \distrfo{M} which denotes the generalized center of the Poisson algebra \smooth{M} in \distrf{M}.
Particularly, we have that $\mathcal{E}(1)\in\distrfo{M}$
$$
\scalar{\brck{f}{\mathcal{E}(1)}}{\alpha}=-\int\limits_M\lder{X_f}{\alpha}=-\int\limits_Md(i_{X_f}\alpha)=0,\quad\all f\in\smooth{M}
$$
\begin{prop}
If the Poisson structure on $M$ is nondegenerate then \distrfo{M} is 1-dimensional an is spanned by $\mathcal{E}(1)$.
\end{prop}
\begin{proof}
As the Poisson structure is nondegenerate, it is given with a symplectic form $\omega$. If $\Phi\in\distrfo{M}$, then
we have the following
$$
\begin{array}{c}
\brck{f}{\Phi}=0,\quad\all f\in\smooth{M}\gamomd\scalar{\brck{f}{\Phi}}{\alpha}=\scalar{\brck{f}{\Phi}}{\psi\cdot\omega^{n/2}}=\\
\\
=-\scalar{\Phi}{\lder{X_f}{\psi\cdot\omega^{n/2}}}=-\scalar{\Phi}{\brck{f}{\psi}\cdot\omega^{n/2}}=0,\quad\all\psi,f\in\smooth{M}
\end{array}
$$
where $n=\dim(M)$. Hence we obtain that the elements of \distrfo{M} vanish on the subspace
$\brck{\smooth{M}}{\smooth{M}}\cdot\omega^{n/2}$ of \forms{n}{M}. But we have that $\forms{n}{M}\cong\smooth{M}\cdot\omega^{n/2}$
and therefore \distrfo{M} is the dual space of
$\frac{\smooth{M}\cdot\omega^{n/2}}{\brck{\smooth{M}}{\smooth{M}}\cdot\omega^{n/2}}$, which itself is isomorphic to
$\frac{\smooth{M}}{\brck{\smooth{M}}{\smooth{M}}}$. As it was mentioned early, in the case of a symplectic manifold,
the latter quotient is isomorphic to $\Real$.
\end{proof}
When the Poisson structure is degenerate, unlike the case of \distr{M}, for \distrf{M} we have no possibilities to
construct distributions like $\delta_N$. But we shall describe another method of constructing of some generalized
Casimirs in \distrfo{M}, besides the constants: $k\cdot\mathcal{E}(1),\ k\in\Real$.

As in the case of \distr{M}, choose a compact symplectic leaf $N$ in $M$. Let $u\in\Gamma(\wedge^{n-2k}V(N))$ be a transversal
multivector field (see section \ref{bott-for-trans-multi}) of degree $n-2k$, where $2k=\dim(N),\ n=\dim(M)$. Consider
the functional on \forms{n}{M} defined by the formula
$$
\scalar{\delta^u}{\alpha}=\int\limits_Ni_u\alpha,\quad\alpha\in\forms{n}{M}
$$
The natural meaning of the inner product $i_u\alpha$ is clear, and after this operation we obtain a $2k$-form on
$N$. The bracket of a function $f\in\smooth{M}$ and the distribution $\delta^u$ gives
$$
\scalar{\brck{f}{\delta^u}}{\alpha}=-\scalar{\delta^u}{\lder{X_f}{\alpha}}=-\int\limits_Ni_u\lder{X_f}{\alpha}
$$
using the formula \ref{lieschout} for the relation between the Schouten-Nijenhuis bracket and Lie derivation
and the formula \ref{BottSchouten} for the Bott covariant derivation on the bundle of transversal multivectors,
the latter expression gives
$$
\scalar{\brck{f}{\delta^u}}{\alpha}=\int\limits_Ni_{\nabla_{_{df}}(u)}(\omega)-\int\limits_N\underbrace{\lder{X_f}{i_u\omega}}_{d(i_{X_f}i_u\omega)}=
\int\limits_Ni_{\nabla_{_{df}}(u)}(\omega)
$$
where $\nabla_{_{df}}(u)$ is the extension of the Bott covariant derivation on $V(N)$ to the exterior degree
$\wedge^{n-2k}V(N)$, defined in Section \ref{bott-for-trans-multi}.

The last formula implies the following
\begin{prop}
For a given transversal multivector field $u\in\Gamma(\wedge^{n-2k}V(N))$ on a compact symplectic leaf $N$, the
distribution $\delta^u\in\distrf{M}$ is an element of the generalized center if and only if the covariant derivative of $u$ with
respect to the Bott connection on $\wedge^{n-2k}V(N)$ is $0$ (in other words: $u$ is \textbf{flat} with respect to
the Bott connection on $\wedge^{n-2k}V(N)$).
\end{prop}
It is clear that the vector bundle $\wedge^{n-2k}V(N)$ on $N$ is 1-dimensional and therefore, if the equation
$\nabla(u)=0$ has at least one nontrivial solution, then the Bott connection on on $\wedge^{n-2k}V(N)$ is flat.
As we see from the Example \ref{example2} in the section \ref{bott-for-trans-multi}, the existence of such section,
in general, is problematic, even in the case when $N$ consists of just one point. And if such section exists, then
it gives only one additional dimension for the space \distrfo{M}, besides the constants. Thus it seems that in this
realization of the test-objects space, at least one dimension for the generalized center is ''lost``. Let us consider
this situation in more details
\begin{exmp}
Let $(S,\ \brck{\ }{\ }_1=\varphi\cdot\brck{\ }{\ })$ be the same as in the example \ref{example}. It follows from
the definition that
$$
\distrfo{S}=\set{\Phi\in\distrf{M}\ |\ \Phi\pars{\brck{f}{\alpha}}=0,\ \all f\in\smooth{M},\ \all\alpha\in\forms{2}{S}}
$$
Otherwise, \distrfo{M} is the orthogonal space of the subspace
$$
\set{\lder{X^1_f}{\alpha}\ |\ f\in\smooth{S},\ \alpha\in\forms{2}{S}}\subset\forms{2}{S}
$$
Any $\alpha\in\forms{2}{S}$ can be represented as $\alpha=\psi\cdot\omega$, where $\omega$ is the symplectic form
corresponding to the original Poisson bracket (which is nondegenerate). For $\alpha$ and $f\in\smooth{S}$ we have
$$
\begin{array}{c}
\lder{X^1_f}{\alpha}=\lder{X^1_f}{\psi\cdot\omega}=\brck{f}{\psi}_1\cdot\omega+\psi\lder{X^1_f}{\omega}=\\
\\
=\varphi\brck{f}{\psi}\cdot\omega+\psi\cdot\lder{{\varphi\cdot X_f}}{\omega}=\varphi\brck{f}{\psi}\cdot\omega+\psi\cdot d\varphi\wedge(i_{X_f}\omega)=\\
\\
=\pars{\varphi\brck{f}{\psi}+\psi\brck{f}{\varphi}}\cdot\omega=\brck{f}{\psi\varphi}\cdot\omega\\
\\
\textrm{(we used the equality }i_{X_f}(\underbrace{d\varphi\wedge\omega}_0)=\brck{f}{\varphi}\cdot\omega-d\varphi\wedge(i_{X_f}\omega))
\end{array}
$$
thus we obtain that \distrfo{S} is the dual to quotient space
$$
\frac{\smooth{S}\cdot\omega}{\brck{\smooth{S}}{\ \varphi\cdot\smooth{S}}\cdot\omega}\cong
\frac{\smooth{S}}{\brck{\smooth{S}}{\ \varphi\cdot\smooth{S}}}
$$
If the point $x_0=\varphi^{-1}(0)$ is not singular for $\varphi$ ($\varphi'_{x_0}\neq0$) then any $f\in\smooth{S}$
can be represented as $f=\varphi\cdot g+f(x_0)$ for some $g\in\smooth{S}$. Therefore in this case we have
$$
\brck{\smooth{S}}{\varphi\cdot\smooth{S}}=\brck{\smooth{S}}{I_{x_0}}
$$
where $I_{x_0}$ is the ideal of functions vanishing at $x_0$. But on the other hand
$$
\brck{\smooth{S}}{I_{x_0}}=\brck{\smooth{S}}{\smooth{S}}
$$
because $\brck{\smooth{S}}{\Real}=\set{0}$. As the original bracket $\brck{\ }{\ }$ is nondegenerate, we have that
$\distrfo{S}\cong\Real$, which explains the ''losing`` one dimension in \distrfo{M} in this realization of the
space of test-objects for distributions.
\end{exmp}
% ----------------------------------------------------------------

\end{document}